\documentclass[11pt]{amsart}
\usepackage{amssymb,amsmath,amsthm}
\usepackage{a4}
\usepackage{epsf}
\usepackage{epsfig}
\usepackage{url}
\begin{document} 

\overfullrule=3pt

\newcommand{\NP}{$\mathcal{N}\mathcal{P}$}
\newcommand{\newsymb}{{\mathcal P}}
\newcommand{\Pn}{{\mathcal P}^n}
\newcommand{\R}{{\mathbb R}}
\newcommand{\N}{{\mathbb N}}
\newcommand{\Q}{{\mathbb Q}}
\newcommand{\Z}{{\mathbb Z}}
\newcommand{\C}{{\mathbb C}}

\newcommand{\enorm}[1]{\Vert #1\Vert}
\newcommand{\inter}{\mathrm{int}}
\newcommand{\conv}{\mathrm{conv}}
\newcommand{\aff}{\mathrm{aff}}
\newcommand{\lin}{\mathrm{lin}}
\newcommand{\cone}{\mathrm{cone}}
\newcommand{\bd}{\mathrm{bd}}

\newcommand{\dist}{\mathrm{dist}}
\newcommand{\trans}{\intercal}
\newcommand{\diam}{\mathrm{diam}}
\newcommand{\vol}{\mathrm{vol}}
\newcommand{\LE}{\mathrm{G}}
\newcommand{\lE}{\mathrm{g}}

\newcommand{\pp}{\mathfrak{p}}
\newcommand{\pf}{\mathfrak{f}}
\newcommand{\pg}{\mathfrak{g}}
\newcommand{\PP}{\mathfrak{P}}
\newcommand{\pl}{\mathfrak{l}}
\newcommand{\pv}{\mathfrak{v}}
\newcommand{\cl}{\mathrm{cl}}
\newcommand{\bx}{\overline{x}}

\def\ip(#1,#2){#1\cdot#2}

\newtheorem{theorem}{Theorem}[section]
\newtheorem{theorem*}{Theorem}
\newtheorem{corollary}[theorem]{Corollary}
\newtheorem{lemma}[theorem]{Lemma}
\newtheorem{remark}[theorem]{Remark}
\newtheorem{definition}[theorem]{Definition}  
\newtheorem{proposition}[theorem]{Proposition}  
\newtheorem{claim}[theorem]{Claim}
\newtheorem{problem}[theorem]{Problem}
\numberwithin{equation}{section}

\title[]{Ehrhart polynomial and Successive Minima} 
%\date{\today}
%\thanks{}
\author{Martin Henk}
\address{Martin Henk and Achill Sch\"urmann, Universit\"at Magdeburg, IAG,
  Universit\"ats\-platz 2, D-39106 Magdeburg, Germany}
\email{\{henk,achill\}@math.uni-magdeburg.de}
\author{Achill Sch\"urmann}
%\address{Universit\"at Magdeburg, Institut f\"ur Algebra und Geometrie,
%  Universit\"ats\-platz 2, D-39106 Magdeburg, Germany}
%\email{achill@math.uni-magdeburg.de}
\author{J\"org M.~Wills}
\address{J\"org M.~Wills, Universit\"at Siegen,  Mathematisches Institut,
 ENC, D-57068 Siegen, Germany}
\email{wills@mathematik.uni-siegen.de}
\dedicatory{Dedicated to Rolf Schneider on the occasion of his 65th birthday}

%\keywords{Lattice polytopes, succesive minima, Ehrhart polynomial}
%\subjclass[2000]{52C07, 11H06}

\begin{abstract} We investigate the Ehrhart polynomial for the class of $0$-sym\-metric convex lattice polytopes in Euclidean $n$-space $\R^n$. It turns out that the roots of the Ehrhart polynomial and Minkowski's successive minima are closely related by their geometric and arithmetic mean.   
We also  show that the roots of lattice $n$-polytopes with or without interior lattice points differ essentially. Furthermore, we study the structure of the roots in the planar case. Here it turns out that their distribution reflects basic properties of lattice polygons.     
\end{abstract}

\maketitle

\section{Introduction}
Let $\mathcal{P}^n$ be the set of all convex lattice $n$-polytopes in the $n$-dimensional Euclidean space $\R^n$ with respect to the standard lattice $\Z^n$, i.e., all vertices of $P\in \mathcal{P}^n$ have  integral coordinates and $\dim(P)=n$. The lattice point enumerator of a set $S\subset\R^n$ is denoted by $\LE(S)$,  i.e., $\LE(S)=\#(S\cap \Z^n)$. Here we restrict ourselves to the lattice $\Z^n$, but  all definitions and results can easily be generalised to arbitrary lattices.

In 1962 Ehrhart \cite{ehrhart:polynomial} showed that for $k\in \N$ the lattice point enumerator $\LE(k\,P)$, $P\in \mathcal{P}^n$, is a polynomial of degree $n$ in $k$ where the coefficients $\LE_i(P)$, $0\leq i\leq n$,  depend only on $P$:
\begin{equation*}
   \LE(k\,P)=\sum_{i=0}^n \LE_i(P)\,k^i.
\end{equation*} 
Moreover in \cite{ehrhart:reciprocity} he  proved his famous ``reciprocity law''  which says that 
\begin{equation}
   \LE(\inter(k\,P))=(-1)^{n} \sum_{i=0}^n \LE_i(P)\,(-k)^i,
\label{eq:reciprocitylaw}
\end{equation}
where $\inter()$ denotes the interior. Two  of the $n+1$ coefficients $\LE_i(P)$ are obvious, namely, $\LE_0(P)=1$ and $\LE_n(P)=\vol(P)$, where $\vol()$ denotes the volume, i.e., the $n$-dimensional Lebesgue measure on $\R^n$.  Also the second leading coefficient admits a simple geometric interpretation  which we present in detail in \eqref{eq:2leadterm}. 
All other coefficients $\LE_i(P)$, $1\leq i\leq n-2$, have no such
direct geometric meanings, except for special classes of polytopes
(cf., e.g., \cite{mordell:dedkind, pommersheim:tetrahedra, diazrobins:latticepolytope, betkegritzmann:valuations, beckpixton:birkhoff, liu:cyclic_ehrhart, hibi:92, mustatepayne:bettyehrhart}). 

Since its discovery the Ehrhart polynomial and its coefficients play an essential role in discrete geometry, geometry of numbers and combinatorics (cf., e.g., \cite{gruber-lekkerkerker:87, stanley:enumcombI, hibi:algebraic_combinatorics, GriWil:latticepoints, lagarias:pointlattices, beckrobins:discrete_continously}). For instance, in \cite{betkekneser:latticehadwiger} Betke and Kneser  showed that the coefficients form a basis of all additive and unimodular invariant functionals on the space $\mathcal{P}^n$.   Stanley  studied the Hilbert series of an integral polytope and proved in this context his famous nonnegativity theorem  \cite{stanley:nonnegative}. 
 For representations  of $\LE_i(P)$ in terms of Todd classes of a toric variety associated with $P$ we refer to \cite{barvinokpommersheim:algorithmic} and the 
references within. Based on Barvinok's methods for counting  lattice points (cf., e.g., \cite{barvinok:ehrhart, barvinok:polynomialtime}),  
  De Loera et al.  developed an efficient algorithm for calculating  the coefficients of the Ehrhart polynomial (cf.~\cite{latte:program, deloera:effecticve_counting}). There were also some attempts to bound $\LE(P)$ in terms of the intrinsic volumes (see \cite{GriWil:latticepoints, wills:90a} as general references), but only with limited success. For relations among the coefficients of the Ehrhart polynomial we refer to \cite{betke-mcMullen:85, becketall:roots_ehrhart}.  

In recent years the Ehrhart polynomial was not only regarded as a polynomial for integers $k$, but as a formal polynomial of a complex variable $s\in \C$ (cf.~\cite{wills:90c, becketall:roots_ehrhart, rodriguez:zeros}). Therefore, for $P\in\Pn$ and $s\in \C$ we set   
\begin{equation*}
   \LE(s, P)= \sum_{i=0}^n \LE_i(P)\, s^i = \prod_{i=1}^n \left(1+\frac{s}{\gamma_i(P)}\right), 
\end{equation*}
where $-\gamma_i(P)\in \C$, $1\leq i\leq n$,  are the roots of the Ehrhart polynomial $\LE(s,P)$. So $\LE_i(P)$ is the $i$-th elementary function of $1/\gamma_1,\dots,1/\gamma_n$, and, in particular, we have 
\begin{equation}
  \vol(P)=\LE_n(P)=\prod_{i=1}^n \frac{1}{\gamma_i(P)}\quad\text{and}\quad 
  \LE_{n-1}(P)=\sum_{j=1}^n \prod_{i\ne j} \frac{1}{\gamma_i(P)}.
\label{eq:coeffroots}
\end{equation}

Here we are interested in relations between $\gamma_i(P)$ and Minkowski's successive minima. To this end let $\mathcal{K}_0^n$ be the class of $0$-symmetric convex bodies in $\R^n$ having non-empty interior.  For $K\in \mathcal{K}_0^n$ and $1\leq i\leq n$  the $i$-th successive minimum $\lambda_i(K)$ is defined by \cite[pp.~58]{gruber-lekkerkerker:87}
\begin{equation*}
    \lambda_i(K)=\min\left\{\lambda>0 : \dim(\lambda\,K\cap\Z^n)\geq i\right\}. 
\end{equation*}

Minkowski proved a  lower and upper bound on the volume of $K$ in terms  of the successive minima which can be formulated as (cf.~\cite[pp.~199]{minkowski:96}, \cite[pp.~58]{gruber-lekkerkerker:87})
\begin{theorem}[Minkowski, 1896] Let $K\in\mathcal{K}_0^n$. Then 
\begin{equation}
     \frac{1}{n!} \prod_{i=1}^n \frac{2}{\lambda_i(K)} \leq \vol(K)\leq \prod_{i=1}^n \frac{2}{\lambda_i(K)}
\label{eq:minkowski}
\end{equation}
and both bounds are tight.
\end{theorem}
For instance, the upper bound is attained for the cube $C_n=\{x\in\R^n : |x_i|\leq 1\}$ of edge length $2$, whereas the cross polytope $C_n^*=\{x\in\R^n : \sum|x_i|\leq 1\}$ matches the lower bound. 
 In contrast to the lower bound the upper bound is a rather deep result in geometry of numbers and has fascinated many mathematicians (cf., e.g.,~\cite{bwz:65, davenport:39, siegel:89, weyl:42}). For extensions of Minkowski's theorems and inequalities between the lattice point enumerator and the successive minima see, e.g., \cite{henk:90, bhw:93, henk:02}.

In terms of the roots of the Ehrhart polynomial (see \eqref{eq:coeffroots}) we can rewrite \eqref{eq:minkowski} as 
\begin{equation}
     \left(\prod_{i=1}^n \frac{\lambda_i(P)}{2}\right)^{1/n} \leq \left(\prod_{i=1}^n \gamma_i(P)\right)^{1/n}\leq  n!^{1/n}\,\left(\prod_{i=1}^n \frac{\lambda_i(P)}{2}\right)^{1/n}
\label{eq:minkroot}
\end{equation}
for the class of all $0$-symmetric lattice $n$-polytopes which we denote by $\mathcal{P}_0^n$.

These inequalities between the geometric mean of the successive minima and the roots of the Ehrhart polynomial lead naturally to the question of further inequalities. Here our main result describes a relation between the arithmetic mean of $\lambda_i(P)/2$ and $\gamma_i(P)$. 
\begin{theorem} Let $P\in \mathcal{P}_0^n$. Then 
\begin{equation}
     \frac{1}{n}\left(\sum_{i=1}^n \gamma_i(P)\right) \leq \frac{1}{n}\left(\sum_{i=1}^n \frac{\lambda_i(P)}{2}\right) 
\label{eq:main}
\end{equation}
and the bound is tight.
\label{thm:main}
\end{theorem}
An interesting feature of this inequality is the fact that it is tight for the cube $C_n$ as well as for the cross polytope $C_n^*$.  
Comparing the equality case in \eqref{eq:main} and   in Minkowski's inequalities  \eqref{eq:minkroot} we found that all possible cases can occur.
\begin{proposition} \hfill \renewcommand{\labelenumi}{\rm \roman{enumi})}
\begin{enumerate}
\item  There exist lattice polytopes $P\in\mathcal{P}_0^n$ with equality in \eqref{eq:main} and equality in the upper bound or in the lower bound or in none of the inequalities of \eqref{eq:minkroot}.   
\item There exist lattice polytopes $P\in\mathcal{P}_0^n$ with strict inequality in \eqref{eq:main} and equality in the upper bound or in the lower bound or in none of the inequalities of \eqref{eq:minkroot}.
\end{enumerate}   
\label{prop:classeq}
\end{proposition}  
Since the successive minima form an increasing sequence 
Theorem \ref{thm:main}  immediately implies   
\begin{corollary} 
Let $P\in \mathcal{P}_0^n$. Then
\begin{equation*}
  \frac{1}{n}\left(\sum_{i=1}^n \gamma_i(P)\right) \leq \frac{\lambda_n(P)}{2}.
\end{equation*} 
\label{cor:analog_one}
\end{corollary}
This bound may  be regarded as an analogue to Minkowski's first theorem on successive minima (cf.~\cite[pp.~75]{minkowski:96}, \cite[pp.~59]{gruber-lekkerkerker:87}) which states $\vol(K)\leq (2/\lambda_1(K))^n$ and which can be rewritten for $P\in \mathcal{P}_0^n$  as 
\begin{equation*}
  \frac{\lambda_1(P)}{2} \leq \left(\prod_{i=1}^n \gamma_i(P)\right)^{1/n}.
\label{eq:analog_one}
\end{equation*} 
Another consequence of Theorem is the following upper bound on  $\LE_{n-1}(P)$

\begin{corollary}  Let $P\in \mathcal{P}_0^n$. Then 
\begin{equation*}
  \LE_{n-1}(P) \leq \sum_{j=1}^n \prod_{i\ne j} \frac{2}{\lambda_i(P)}
\end{equation*}
and the bound is tight.
\label{cor:secondcoeff}
\end{corollary}
To see this we note that  
in terms of the coefficient $\LE_{n-1}()$ inequality \eqref{eq:main} is equivalent to (cf.~\eqref{eq:coeffroots})
\begin{equation}
     \frac{\LE_{n-1}(P)}{\vol(P)} \leq \sum_{i=1}^n \frac{\lambda_i(P)}{2}.
\label{eq:whattoshow}      
\end{equation}  
Hence, together with the upper bound  in  \eqref{eq:minkowski}  we get Corollary \ref{cor:secondcoeff}, where equality is attained, e.g., for the cube $C_n$. For the special case $\lambda_i(P)=1$ Corollary \ref{cor:secondcoeff} gives $ \LE_{n-1}(P) \leq n\,2^{n-1}$ which was already shown in   \cite{wills:minkowskianalog}.

In order to study the size and the distribution of the zeros of Ehrhart polynomials we introduce 
\begin{definition} \hfill \renewcommand{\labelenumi}{\rm \roman{enumi})}
\begin{enumerate} 
\item For $l\in\N$ let $\mathcal{P}^n(l)$ ($\mathcal{P}^n_0(l)$) be the set of ($0$-symmetric) lattice polytopes $P$ having exactly $l$ interior lattice points, i.e., $\LE(\inter(P))=l$.
\item For $l\in\N$ let $\Gamma(n,l)$ ($\Gamma_0(n,l)$) be the zeros  of Ehrhart polynomials with respect to   $\mathcal{P}^n(l)$ ($\mathcal{P}^n_0(l)$).
\item Let $\Gamma(n)=\cup_{l=0}^\infty \Gamma(n,l)$ ($\Gamma_0(n)=\cup_{l=0}^\infty \Gamma_0(n,l)$) be the set of all
  zeros of Ehrhart polynomials w.r.t.~ $n$-dimensional ($0$-symmetric)
  lattice polytopes.  
\end{enumerate}
\end{definition}

In \cite[Theorem 1.2 (a)]{becketall:roots_ehrhart} the general upper bound $|s|\leq (n+1)!+1$, $s\in\Gamma(n)$,
was proven. This result reflects the fact that lattice polytopes
cannot be too small.  Here we consider the case that  lattice
polytopes cannot be too large by fixing the number $l\in\N$ of interior lattice
points. It turns out that the case $l=0$ is completely different from $l\geq 1$.
\begin{theorem} \hfill \renewcommand{\labelenumi}{\rm \roman{enumi})}
\begin{enumerate}
\item $\Gamma(n-1,0)\subset \Gamma(n,0)$ and $0$ is a cluster point of $\Gamma(n,0)$. Furthermore,  $1$ is also a cluster point of $\Gamma(n,0)$ for $n\geq 3$. 
In particular, the set $\Gamma(n,0)$ is infinite. 
\item $\Gamma(n-1,l)\subset \Gamma(n,l)$ and $\Gamma(n,l)$ is finite for any $l\geq 1$.
%\item For each $n$ and $k\in\N$ there exist polytopes $P_j\in \mathcal{P}^n(k)$, $j\in\N$, such that $\inf\{|\gamma_i(P_j)| : 1\leq i\leq n, \, j\in\N\} >0$.
\item For $P\in  \mathcal{P}^n_0(l)$ we have 
\begin{equation*}
   \frac{1}{2} \left(\frac{2}{l+1}\right)^{1/n} \leq \left(\gamma_1(P)\cdot\ldots\cdot\gamma_n(P)\right)^{1/n}
\end{equation*}
and the bound is tight.
\end{enumerate}
\label{thm:boundsinterior}
\end{theorem}

In \cite[Theorem 1.2]{becketall:roots_ehrhart} it was also shown that the real roots
$s\in\Gamma(n)\cap\R$ satisfy 
$s\in[-n,\lfloor n/2\rfloor)$. For $n\leq 4$ the upper bound was improved to
$1$ and  Theorem \ref{thm:boundsinterior} i) implies that this bound is best possible. 
 The lower bound $-n$ for the roots of an Ehrhart polynomial is in general best possible, but it is quite likely that it can
be improved for the class of $0$-symmetric polytopes as the next proposition indicates. 

\begin{proposition} For $n\leq 3$ we have $\gamma\geq -1$ for all $\gamma\in \Gamma_0(n)\cap\R$. 
\label{prop:bound3}
\end{proposition}
Obviously, in the two dimensional
case much more is known and in \cite[Theorem 2.2]{becketall:roots_ehrhart} bounds for the set $\Gamma(2)$ are given. The following result gives information about the strcuture of $\Gamma(2)$.%\newpage
\begin{theorem} \hfill \renewcommand{\labelenumi}{\rm \roman{enumi})}
\begin{enumerate} 
\item The complex roots $\gamma\in \Gamma(2)\cap\C$ lie on the circles
 \begin{equation*}
 \left\{s\in \C : \left|s+2/\LE(\bd P)\right| = 2/\LE(\bd P)\right\}, 
\end{equation*}
where $\bd P$  denotes the boundary. \label{enum:complex}
\item The lattice polygons with real roots are exactly those satisfying 
\begin{equation*}
     \left(\LE(\bd P)/4-1 \right)^2\geq \LE(\inter P).
\end{equation*}
\item 
\begin{equation*}
          \Gamma(2)\subset\{-2,-1,-2/3\}\cup\left\{s=a+\imath b\in \C :
          -\frac{1}{2}\leq a <0, |s+\frac{2}{3}| \leq \frac{2}{3}\right\}.
\end{equation*} \label{enum:all}
\item The cluster points of $\Gamma(2)$ are all  points in $[-1/2,0]$. \label{enum:cluster} 
\item $\Gamma(2,0)=\{ -1, -2/l: l\in\N\}$. \label{enum:cluster0}
\item On the line $\Re s=-1/2$ are just the roots of lattice polygons $P$ with $\LE(\inter P)=1$ and $\LE(\bd P)\leq 8$. The only remaining $P$ with $\LE(\inter(P))=1$ has 9 boundary lattice points and roots $-2/3$, $-1/3$. \label{enum:line} 
%There are exactly 11 roots $a+\imath\, b$ on the line $a = -1/2$.
\end{enumerate}
\label{thm:twodim_roots}
\end{theorem}

Property ii) of the theorem above may be interpreted as an isoperimetric inequality for the lattice point enumerator. 
The following figure depicts the roots (with $\Re s\leq 2/3$) of two dimensional lattice polygons.

\begin{figure}[htb]
\begin{center}
\epsfig{file=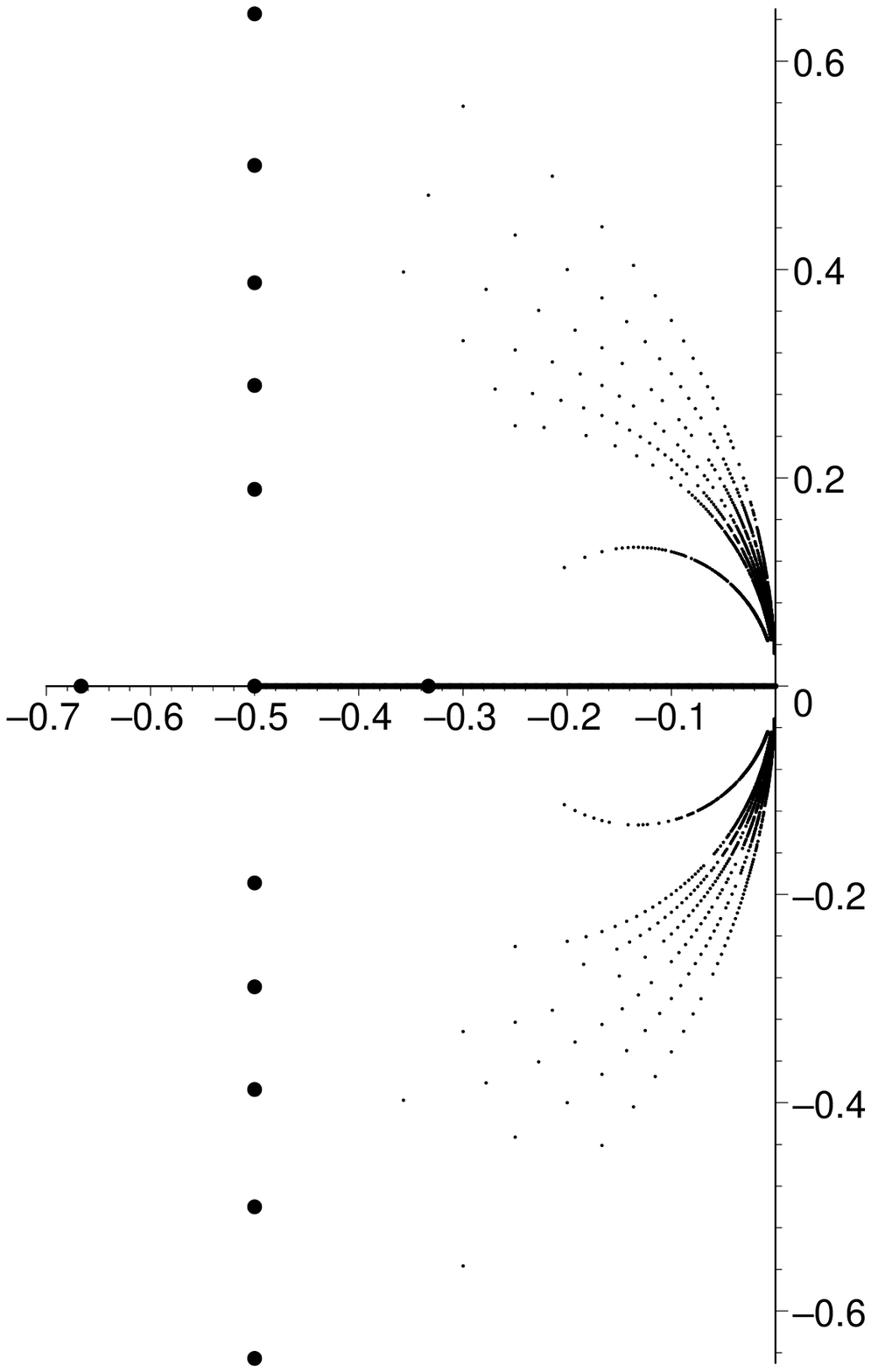, width=4.5cm}
\end{center}
\caption{}
\end{figure}

The paper is organised as follows. In Section 2 we present some basic
facts about the roots of  Ehrhart polynomials and what is known for
some special lattice polytopes. The proofs of Theorem \ref{thm:main}
and Proposition \ref{prop:classeq} are given in Section 3.  Section 4
deals with bounds for the roots, i.e.,
with Theorems \ref{thm:boundsinterior}, \ref{thm:twodim_roots} and Proposition \ref{prop:bound3}. Some open problems and possible generalisations of our
results are given in the last section.
 
\section{Basic facts and examples}
It was already shown by Ehrhart \cite{ehrhart:polynomialII} that 
  the second leading coefficient $\LE_{n-1}(P)$ of the Ehrhart polynomial admits a simple geometric interpretation via the facets  $F_1,\dots,F_m$ of a lattice polytope $P\in\mathcal{P}^n$. Namely, 
\begin{equation}
   \LE_{n-1}(P) = \frac{1}{2} \sum_{i=1}^m \frac{\vol_{n-1}(F_i)}{\det(\aff F_i\cap\Z^n)},
\label{eq:2leadterm}
\end{equation}
where $\vol_{n-1}()$ denotes the $(n-1)$-dimensional volume and 
$\det(\aff F_i\cap\Z^n)$ denotes the determinant of the $(n-1)$-dimensional sublattice of $\Z^n$ contained in the affine hull of the facet $F_i$. So $\LE_{n-1}(P)$ may be regarded as the normalised surface area of $P$ with respect to $\Z^n$.

For instance, for the standard simplex $T_n=\{x\in\R^n : x_i\geq 0,\,\sum x_i\leq 1\}$ we have $\LE(s,T^n)=\binom{n+s}{n}$ and so we get 
\begin{equation*}
  \LE_n(T^n)= \frac{1}{n!}\quad\text{and}\quad \LE_{n-1}(T^{n-1})=\frac{1}{2}\frac{n+1}{(n-1)!}.
\end{equation*}
Since $\LE(\inter(kT_n))=0$ for $k\in\{1,\dots,n\}$ we conclude by Ehrhart's reciprocity law \eqref{eq:reciprocitylaw} that the roots of $\LE(s,T_n)$ are given by the integers $\{-n,\dots,-1\}$. 

Obviously, $\LE(s,P)\ne 0$ for $s\in\N$ and 
if $P$  contains an interior lattice point then its Ehrhart polynomial has no integer roots at all.   
For odd dimensions the polynomial $\LE(s,P)$ has always one real negative root, but it may have positive real roots (see Theorem \ref{thm:boundsinterior} i)). 

For a positive integer vector $m\in\N$ with $m_1\leq m_2\leq \cdots\leq m_n$ let $Q_n(m)=\{x\in\R^n : |x_i| \leq m_i\}$ be the lattice box of edge lengths $2\,m_1,\dots,2\,m_n$. Then 
\begin{equation*}
  \LE(s,Q_n(m)) = \prod_{i=1}^n \left(1+2\,m_i\,s\right) 
\end{equation*}  
and we find 
\begin{equation}
    \frac{\lambda_i(Q_n(m))}{2}=\frac{1}{2\,m_i}=\gamma_i(Q_n(m)),\,1\leq i\leq n.
\label{eq:latticebox}
\end{equation}
So we have a simple relation between the roots of the Ehrhart polynomial and the successive minima. In particular, for the cube $C_n$ we  obtain 
\begin{equation}
\LE_i(C_n)=2^{n-i}\binom{n}{i}.
\label{eq:cube}
\end{equation}
 The coefficients of the Ehrhart polynomial of the cross polytope $C_n^*$ are not so easy to describe. Here we have  
\begin{equation*}
  \LE(s, C_n^*) = \sum_{i=0}^n 2^i \binom{n}{i}\binom{s}{i}
\end{equation*} 
and so 
\begin{equation}
   \LE_n(C_n^*)=\frac{2^n}{n!},\quad\LE_{n-1}(C_n^*)=\frac{2^{n-1}}{(n-1)!}\quad\text{and}\quad \LE_{n-2}(C_n^*)=\frac{2^{n-2}}{(n-2)!}\frac{n+1}{3}.
\label{eq:crosspolytope}
\end{equation}
It follows from a much more general result of Bump et al. \cite[Theorem 4]{bump:crosspolytope} that the zeros of $\LE(s, C_n^*)$ all have real parts equal to $-1/2$. They also prove the functional equation $\LE(-s,C_n^*)=(-1)^n \LE(s-1,C_n^*)$ (cf.~\cite{rodriguez:zeros}). 

Next we note that for two lattice polytopes $P_1\in\mathcal{P}^{n_1}$,
$P_2\in \mathcal{P}^{n_2}$ living in complementary spaces $\R^{n_i}$ we
have for $P_1\oplus P_2\subset  \R^{n_1}\oplus\R^{n_2}$ 
\begin{equation}
 \LE(s,P_1\oplus P_2)=\LE(s,P_1)\,\LE(s,P_2).
\label{eq:directzeros} 
\end{equation} 
In particular the roots of $\LE(s,P_1\oplus P_2)$ are given by the roots of
$\LE(s,P_i)$, $i=1,2$.

Finally we mention that Ehrhart's theorem on the polynomial behaviour of $\LE(k\,P)$ has a well known predecessor in the planar case, namely the so called Pick's theorem \cite{pick:lattice} which states for $P\in\mathcal{P}^2$ and $k\in\N$  (cf.~\eqref{eq:2leadterm})
\begin{equation}
 \LE(k\,P) = \vol(P)\,k^2 + \frac{1}{2}\#(\bd P\cap\Z^2)\,k +1.
\label{eq:pick}
\end{equation}

\section{Successive Minima and the 
         roots  of  the Ehrhart polynomial}
In the following let $P=\{x\in\R^n : a_j\,x\leq b_j,\,1\leq j\leq m\}$ 
be a $0$-symmetric  $n$-dimensional polytope with $a_j\in\R^n$ and $b_j\in\R_{>0}$. The facet corresponding to $a_j$ is denoted by $F_j$ and hence we may write  by the ``pyramid-formula''  
\begin{equation}
     \vol(P)  = \sum_{j=1}^m \frac{\vol_{n-1}(F_j)}{\enorm{a_j}}\frac{b_j}{n}.
\label{eq:vol}
\end{equation}
Here $\enorm{\cdot}$ denotes the Euclidean norm. For the proof of Theorem \ref{thm:main} we need the following lemma. 
\begin{lemma} Let $v_1,\dots,v_n\in\R^n$ be linearly independent. For $1\leq k\leq n-1$ let 
$V_k=\{ j : v_i\,a_j=0,\,1\leq i\leq k\}$ and 
\begin{equation*}
 \nu_k=\sum_{j\in V_k} \frac{\vol_{n-1}(F_j)}{\enorm{a_j}}\frac{b_j}{n}
\end{equation*}
that part of the volume of $P$ corresponding to the pyramids with bases $F_j$, $j\in V_k$. Then we have 
\begin{equation*}
     \vol(P)\geq  \frac{n}{n-k}\nu_k
\end{equation*}
and the bound is tight.
\label{lem:main}
\end{lemma} 
\begin{proof} The cube $C_n$ and the unit vectors $e_i$ as vectors $v_i$ show that the bound is best possible.  Let us fix $k$ such that $V_k\ne\emptyset$ and let $L_k=\lin\{v_1,\dots,v_k\}$ with orthogonal complement $L_k^\perp$. 
Since $P$ is $0$-symmetric we have 
\begin{equation}
     \vol_k(P\cap L_k) \geq \max_{x\in L_k^\perp}\, \vol_k(P\cap (x+L_k)),
\label{eq:volume_trivial}
\end{equation}
where $\vol_k(\cdot)$ denotes the $k$-dimensional volume. 
For $j\in V_k$ we consider $M_j=\conv\{F_j,P\cap  L_k\} \subset P$. By the definition of $V_k$ we have $\inter M_j\cap \inter M_{\overline j}=\emptyset$ for $j\ne \overline{j}$ and hence it suffices to show 
\begin{equation}
    \vol(M_j) \geq \frac{n}{n-k}\,\cdot \frac{\vol_{n-1}(F_j)}{\enorm{a_j}}\frac{b_j}{n} 
\label{eq:toshow}
\end{equation}
for each $j \in V_k$. To this end we may assume $v_i=e_i$, $1\leq i\leq k$, 
$a_j=e_{k+1}$ and $b_j=1$. Then we can write  
\begin{equation*}
\vol(M_j) =\int\limits_{x_{k+1}=0}^1 \int_{x_{k+2}}\dots \int_{x_{n}} \vol_k\left( M_j\cap ( [x_{k+1},\dots,x_n] + L_k)\right) {\rm d}\,x_n\dots{\rm d}\,x_{k+1},
\end{equation*} 
where $[x_{k+1},\dots,x_n]$ denotes the vector $(0,\dots,0,x_{k+1},\dots,x_n)^\intercal\in L_k^\perp$. By \eqref{eq:volume_trivial} and the Brunn-Minkowski  theorem (cf.~\cite[pp.~309]{schneider:93}) we get for $x_{k+1}\in(0,1]$
\begin{equation*}
\begin{split}
  &\vol_k\left( M_j\cap ([x_{k+1},\dots,x_n] + L_k)\right)^\frac{1}{k} \\&\geq
     x_{k+1} \vol_k\left( M_j\cap ([1,\frac{x_{k+2}}{x_{k+1}},\dots,\frac{x_n}{x_{k+1}}] + L_k)\right)^\frac{1}{k}  + (1-x_{k+1})\vol_k(P\cap L_k)^\frac{1}{k}\\ 
&\geq  
\vol_k\left( M_j\cap (\left[1,\frac{x_{k+2}}{x_{k+1}},\dots,\frac{x_n}{x_{k+1}}\right] + L_k)\right)^\frac{1}{k}.
\end{split}
\end{equation*}
Hence we have 
\begin{equation*}
\begin{split}
& \vol(M_j)  \geq \\  &\int\limits_{x_{k+1}=0}^1 \int_{x_{k+2}}\dots \int_{x_{n}}\vol_k\left( M_j\cap ([1,\frac{x_{k+2}}{x_{k+1}},\dots,\frac{x_n}{x_{k+1}}] + L_k)\right) {\rm d}\,x_n\dots{\rm d}\,x_{k+1} 
\\ =&\int\limits_{0}^1 \left(\int\limits_{x_{k+2}}\dots \int\limits_{x_{n}} \vol_k\left( M_j\cap (\left[1,x_{k+2},\dots,x_n\right] + L_k) \right){\rm d}\,x_n\dots{\rm d}\,x_{k+2}\right) t^{n-k-1} {\rm d}\,t \\
&\quad\quad\quad\quad = \vol_{n-1}(F_j) \frac{1}{n-k},
\end{split}
\end{equation*}
which shows \eqref{eq:toshow}.
\end{proof} 

Next we come to the proof of Theorem \ref{thm:main}.
\begin{proof}[Proof of Theorem \ref{thm:main}] Let $P=\{x\in\R^n : a_j\,x\leq b_j,\,1\leq j\leq m\}$ be a  $0$-symmetric lattice polytope. We can assume $a_j\in\Z^n$, $b_j\in\N$, and that the vectors $a_j$ are primitive, i.e., $\conv\{0,a_j\}\cap\Z^n=\{0,a_j\}$. 
As mentioned in the introduction the inequality of the theorem is equivalent to (cf.~\eqref{eq:whattoshow}) 
\begin{equation}
         \frac{\LE_{n-1}(P)}{\vol(P)} \leq \sum_{i=1}^n \frac{\lambda_i(P)}{2}.
\label{eq:whattoshow2}     
\end{equation}
Lattice boxes (cf.~\eqref{eq:latticebox}) as well as the cross polytope (cf.~\eqref{eq:crosspolytope}) show that the inequality is tight.  As before let $F_j$ be the facet of $P$ corresponding to the normal vector $a_j$. Since $a_j$ is primitive we have $\enorm{a_j}=\det(\aff F_j\cap\Z^n)$ and hence we may write (cf.~\eqref{eq:2leadterm})
\begin{equation*}
     \frac{\LE_{n-1}(P)}{\vol(P)} =\frac{1}{2}\sum_{j=1}^m
     \frac{\alpha_j}{\vol(P)} \quad\text{ with } \quad \alpha_j=\frac{\vol_{n-1}(F_j)}{\enorm{a_j}}.
\end{equation*}
Let $\lambda_i$, $1\leq i\leq n$, be the successive minima of $P$ and let   
let $v_1\dots,v_n\in P$ be linearly independent such that
$\lambda_i\,v_i=z_i\in\Z^n$, $1\leq i\leq n$. Next we define the sets $V_k$ and the
volumes $\nu_k$ with respect to these vectors $v_i$ as in the Lemma
\ref{lem:main}. In particular, we may write 
\begin{equation*}
     \nu_k=\sum_{j\in V_k} \alpha_j\frac{b_j}{n}.
\end{equation*}

In addition, we set $V_n=\emptyset$, $V_0=\{1,\dots,m\}$ and $\nu_0=\vol(P)$.
So $V_{k}\subset V_{k-1}$, $1\leq k\leq n$, and  let $q$ be
the smallest number with $V_q =\emptyset$.  
Since $z_i\in\lambda_i\,P$ we have for $1\leq j\leq m$, $1\leq i\leq n$, 
\begin{equation*}
   b_j\geq \frac{1}{\lambda_i}|a_j\,z_i|.
\end{equation*} 
Thus by the definition of the sets $V_k$ we obtain for $1\leq k\leq q$
\begin{equation*}
   b_j\geq \frac{1}{\lambda_k},\quad \text{ for all }j\in V_{k-1}\setminus V_k.
\end{equation*}
Hence we  get 
\begin{equation*}
\begin{split}
 \frac{\LE_{n-1}(P)}{\vol(P)} & = \frac{1}{2}\sum_{k=1}^q
 \frac{\sum_{j\in V_{k-1}\setminus V_k}\alpha_j}{\vol(P)}\\
&\leq \frac{n}{2} \sum_{k=1}^q \lambda_k\,\frac{\sum_{j\in V_{k-1}\setminus V_k}\alpha_j\,\frac{b_j}{n}}{\vol(P)}
=  \frac{n}{2} \sum_{k=1}^q \lambda_k \frac{\nu_{k-1}-\nu_k}{\vol(P)} \\ 
& = \frac{n}{2}\left(\lambda_1 + \sum_{k=1}^{q-1}   \frac{\nu_k}{\vol(P)}(\lambda_{k+1}-\lambda_k)\right) .
\end{split}
\end{equation*}
Since the successive minima form an increasing sequence we may apply Lemma \ref{lem:main} and get 
\begin{equation*}
\begin{split}
\frac{\LE_{n-1}(P)}{\vol(P)} & \leq \frac{n}{2}\left(\lambda_1 + \sum_{k=1}^{q-1}   \frac{n-k}{n}(\lambda_{k+1}-\lambda_k)\right)
 \\
& =\frac{1}{2}\left(\sum_{i=1}^{q-1}\lambda_i +(n-q+1)\lambda_q\right) \leq \frac{1}{2}\left( \lambda_1 + \lambda_2 + \cdots + \lambda_n\right).
\end{split}
\end{equation*}
\end{proof}

\begin{remark}{\em  \hfill \renewcommand{\labelenumi}{\rm \roman{enumi})}
\begin{enumerate}
 \item In contrast to the geometric mean, the arithmetic mean of  $\gamma_i(P)$ is not bounded from below by the arithmetic mean of  $2/\lambda_i(P)$. Or to phrase it in another way: 
there is no lower bound on $G_{n-1}(P)/\vol(P)$ in terms of the sum of the successive minima. To see this we consider for an integer $l\in\N$ the so called Kleetope (see \cite[pp.~217]{gruenbaum:polytopes}) 
\begin{equation*}
     M_l = \conv\{l\,C_n, \pm (l+1)\,e_i, 1\leq i\leq n\},
\label{eq:example_nolower}
\end{equation*}
i.e., on each facet of the cube  $l\,C_n$ we place a little pyramid.
Then $\vol(M_l)=(2\,l)^n(1+1/l)$, $\lambda_i(M_l)=1/(l+1)$, $1\leq i\leq n$, and 
based on \eqref{eq:2leadterm} one can calculate $\LE_{n-1}(M_l)= 2\,n\,(2\,l)^{n-2}$.
Hence  
\begin{equation*}
    \frac{\LE_{n-1}(M_l)}{\vol(M_l)}=\frac{1}{l}\sum_{i=1}^n \frac{\lambda_i(M_l)}{2}.
\end{equation*}
\item Looking at the arithmetic and geometric mean one may also ask for further relations between the elementary functions of  $\gamma_i$ and $\lambda_i/2$. There does not seem, however,  to be an obvious relation. For instance for the cross polytope $C_n^\star$ with roots $\gamma_i$ and successive minima $\lambda_i$ we have (cf.~\eqref{eq:coeffroots}, \eqref{eq:crosspolytope})
\begin{equation*}
\sum_{i\ne j} \gamma_i\gamma_j 
=\frac{\LE_{n-2}(C_n^\star)}{\vol(C_n^\star)} 
 = \frac{n+1}{6} \binom{n}{2}=  
 2\frac{n+1}{3}\sum_{i\ne j} \frac{\lambda_i}{2}\frac{\lambda_j}{2}.
\end{equation*}
   \end{enumerate}  
}  
\label{rem:elementary}
\end{remark}

Now we come to Proposition \ref{prop:classeq}, i.e., to the discussion of
the equality cases in Theorem \ref{thm:main} in comparison with those
in Minkowski's inequalities \eqref{eq:minkowski} and \eqref{eq:minkroot}, respectively. 

\begin{proof}[Proof of Proposition \ref{prop:classeq}] 
As mentioned in the proof
above we have equality in \eqref{eq:whattoshow} (and thus in \eqref{eq:main}) for lattices
boxes and crosspolytopes which are also extremal bodies for the  left 
hand and right hand inequality in    \eqref{eq:minkroot}, respectively. However,
there are also bodies with equality in \eqref{eq:main} but strict inequality
in both inequalities of \eqref{eq:minkroot}. To this end we firstly consider in
the planar case the hexagon $H=\conv\{\pm\,e_1,\pm\,e_2,\pm\,(1,1)^\intercal\}$. 

Then we have $\lambda_i(H)=1$, $\vol(H)=3$ and strict inequality
in \eqref{eq:minkroot}. On the other hand  
its Ehrhart polynomial is given by $\LE(s,H)=3\,s^2+3\,s+1$. Hence
$\gamma_1(H)+\gamma_2(H)=1$ and again we have equality in
\eqref{eq:main}. For $n\geq 3$ we consider $P=H\oplus C_{n-2}\in\mathcal{P}^n_0$, i.e., the Cartesian product of $H$ and the
cube $C_{n-2}$. Then we have  $\lambda_i(P)=1$, $1\leq i\leq n$,
$\vol(P)=3\,2^{n-2}$, and since the roots of $P$ are just the roots of
$H$ plus the roots of $C_{n-2}$ (see \eqref{eq:directzeros}), 
the polytope $P$ gives the desired example. Hence we
have shown the first statement of Proposition \ref{prop:classeq}. 

For the second
part we start with a space-filling planar hexagon $\tilde{H}$ with vertices $\{\pm (3,-3)^\intercal, \allowbreak \pm (3,5)^\intercal, \pm (5,3)^\intercal\}$. Then we have $\lambda_i(\tilde{H})=1/4$, $\vol(\tilde{H})=64$ and, in particular, equality in the left hand side of \eqref{eq:minkroot}. On the other hand $ \tilde{H}$ has 12 lattice points on the boundary and so the roots of its Ehrhart polynomial sum up to $-6/64$ (cf.~\eqref{eq:pick}) which shows that we have strict inequality in \eqref{eq:main}. As before we can easily generalise this example to higher dimensions by taking the Cartesian product $\tilde{H}\oplus C_{n-2}$. 
In order to have equality in the right hand side of \eqref{eq:minkroot} we consider the polygon $K=\conv\{\pm e_1, \pm 2\,e_2\}$ of volume $4$. Here we have $\lambda_1(K)=1/2$, $\lambda_2(K)=1$ and equality in the right hand side of  \eqref{eq:minkroot}, but again inequality in \eqref{eq:main}.
For a ``generic'' polytope we have inequality in all three inequalities. As an example we can take the polytope $M_l$ considered in Remark \ref{rem:elementary} i). 
\end{proof}

\section{Bounds for the roots of the Ehrhart polynomial}
% For the proof of Theorem \ref{} we need the so called elongated
% crosspolytopes which for an odd integer $k$ are defined by 
% \begin{equation*}
% C_n^*(k)=\conv\{\pm\,\frac{1}{2}(k+1)e_1,\pm\,e_2,\dots,\pm\,e_n
% \} \in\mathcal{P}^n_0(k).
% \end{equation*}
% Moreover, for an odd $k$ let 
% \begin{equation*}
% Q_n(k)=\{x\in\R^n : |x_1|\leq \frac{1}{2}(k+1),\,|x_i|\leq 1,\, 2\leq i\leq n\} \in\mathcal{P}^n_0(k)
% \end{equation*}
% be the lattice box of edge lengths $k+1$, $2,\dots,2$. 
% Obviously, we have $\vol(C_n^*(k))=(k+1)\,2^{n-1}/n!$ and
% $\vol(Q_n(k))=(k+1)2^{n-1}$.
The main message of Theorem \ref{thm:boundsinterior} is that the roots of Ehrhart polynomials of lattice polytopes $P$ behave differently depending on whether the polytope contains interior lattice points or not.   
\begin{proof}[Proof of Theorem \ref{thm:boundsinterior}] 
First we observe that given an $(n-1)$-dimensional  lattice polytope $P\in\mathcal{P}^{n-1}(l)$,
$l\in\N$, we may embed  $P$ into $\R^n$ such that $P$ lies in
the hyperplane $\{x\in\R^n : x_n=0\}$. Thus the prism
$\overline{P}=P+\conv\{-e_n,e_n\}$ of height $2$ over the basis $P$
belongs to $\mathcal{P}^{n}(l)$ and the roots of $\overline{P}$ are
the roots of $P$ plus $-1/2$ (cf.~\eqref{eq:directzeros}). 
So  we have $\Gamma(n-1,l)\subset\Gamma(n,l)$ for $l\in\N$. 

Hence Theorem \ref{thm:twodim_roots} v)  shows already that $0$ is a cluster point of $\Gamma(n,0)$ and 
in order to show that $\gamma=1$ is a cluster point of $\Gamma(n,0)$ for $n\geq 3$, it suffices to prove it for  $n=3$. Here we consider  for an integer $q$ the so called  Reeve simplex $R_q=\conv\{0,e_1,e_2, e_1+e_2+q\,e_3\}$ (see \cite{reeve:simplex}).
It is easy to calculate that 
\begin{equation*}
  \LE(s,R_q)= \frac{q}{6}\,s^3 + s^2 + \frac{12-q}{6}\,s +1 = \frac{q}{6}\,(s+1)\left(s^2-s\left(1-\frac{6}{q}\right)+\frac{6}{q}\right).
\end{equation*}     
Thus the smallest root is $-1$ and for $q\to\infty$ the two other roots converge to $1$ and $0$, respectively. 
 
For $l\geq 1$ we know that the volume of $P\in \mathcal{P}^n(l)$ is bounded (cf.~\cite{lagariasziegler:volume}) by a constant depending only on $l$ and $n$. Thus, up to unimodular transformations, there are only finitely many different polytopes $P\in \mathcal{P}^n(l)$. Hence $\Gamma(n,l)$ is finite. 
% For an odd integer $l$ let
% \begin{equation*}
%  Q_n(l)=\{x\in\R^n : |x_1|\leq \frac{1}{2}(l+1),\,|x_i|\leq 1,\, 2\leq i\leq n\} \in\mathcal{P}^n_0(l)
% \end{equation*}
% be the lattice box of edge lengths $l+1$, $2,\dots,2$. 
% Obviously, we have $\vol(Q_n(l))=(l+1)2^{n-1}$ and the roots of its Ehrhart polynomial are $-1/(l+1)$ and $n-1$ times $-1/2$.  Thus $Q_{n}(l)\oplus\conv\{0,e_{n+1}\}$ has the same zeros plus $-1$ (see \eqref{eq:directzeros}). Since the latter polytope belongs to $\mathcal{P}^{n+1}(0)$ we have shown that $0$ is  a cluster point of $\Gamma(n,0)$, $n\geq 2$. 
For the third statement we use a result of Blichfeldt  and van der Corput \cite[p.~51]{gruber-lekkerkerker:87} which says that for $P\in \mathcal{P}^n_0(l)$ 
\begin{equation*}
 \vol(P)\leq 2^n\,\frac{l+1}{2}. 
\end{equation*} 
Together with \eqref{eq:coeffroots} we get the lower bound in Theorem \ref{thm:boundsinterior} iii) and the box 
\begin{equation*}
 Q_n(l)=\{x\in\R^n : |x_1|\leq \frac{1}{2}(l+1),\,|x_i|\leq 1,\, 2\leq i\leq n\} \in\mathcal{P}^n_0(l),\quad l \text{ odd},
\end{equation*}
shows that the bound is tight. 
\end{proof}

%\begin{remark}{\em % \hfill \renewcommand{\labelenumi}{\rm \roman{enumi})}
% \begin{enumerate}
% \item 
The last two statements of Theorem \ref{thm:boundsinterior} imply that there are constants $\alpha(n,l),\beta(n,l)>0$, depending only on $l$ and $n$,  such that for $\gamma\in\Gamma(n,l)$ 
\begin{equation*}
          \alpha(n,l)\leq |\gamma|\leq\beta(n,l).
\end{equation*}  
It seems to be, however, rather difficult to give good bounds for these constants. 
% \item 
% \end{enumerate}
%}
%\end{remark}
Furthermore, we remark that numerical computations up to dimension 8 indicate that 1 is  also a cluster point of the roots of higher dimensional Reeve simplices given by $\conv\{0,e_1,\dots,\allowbreak e_{n-1}, \sum_{i=1}^{n-1}e_i+q\,e_n\}$.

\begin{proof}[Proof of Proposition \ref{prop:bound3}] On account of Theorem \ref{thm:boundsinterior} i) and ii) it suffices to prove the statement for $n=3$. 
Let $P\in\mathcal{P}_0^3$ and for abbreviation we write   
$\LE(s)=\LE(s,P)=\lE_3\,s^3+\lE_2\,s^2+\lE_1\,s+1$.  By Ehrhart's reciprocity law \eqref{eq:reciprocitylaw} we  have $-\lE_3+\lE_2-\lE_1+1\leq -1$ or equivalently $\lE_3-\lE_2+\lE_1\geq 2$. Hence computing the derivative yields  
\begin{equation*}
  \LE^\prime(-1) = 3\,\lE_3 - 2\,\lE_2 + \lE_1 \geq 2\,\lE_3-\lE_2 +2 > 2,
\end{equation*}   
where the last inequality follows from \eqref{eq:whattoshow} which in the 3-dimensional case becomes $3\,\lE_3-2\,\lE_2\geq 0$. Moreover we find for $s<-1$ 
\begin{equation*}
  \LE^\prime(s) -\LE^\prime(-1) = 3\,\lE_3(s^2-1)+2\,\lE_2(s+1)=(s+1)\left(3\,\lE_3\,s - (3\,\lE_3-2\,\lE_2)\right) > 0.
\end{equation*}
Thus for $s<-1$ we have $\LE^\prime(s) > \LE^\prime(-1)>2$ and together with $\LE(-1)\leq -1$ we get $\LE(s) \leq -1 $ for $s\leq -1$. Hence all roots are greater than $-1$.
\end{proof}

Now we come to the roots of lattice polygons. 

\begin{proof}[Proof of Theorem \ref{thm:twodim_roots}] 
For part i) let $P$ be a lattice polygon with Ehrhart polynomial $\lE_2\,s^2+\lE_1\,s+1$ and complex roots 
\begin{equation}
  -\gamma_{1,2} = - \frac{\lE_1}{2\,\lE_2}\pm \imath \sqrt{\frac{1}{\lE_2}-\left(\frac{\lE_1}{2\,\lE_2}\right)^2}.
\label{eq:tworoots}
\end{equation}
Since $\lE_1$ is just one half the number of lattice points on the boundary of $P$ (see \eqref{eq:pick})  we set $\lE_1=m/2$, $m\geq 3$. Setting $a=-m/(4\,\lE_2)$ and $b=\sqrt{-4\,a/m-a^2}$ we may write $ -\gamma_{1,2} = a\pm \imath b$ and find  
\begin{equation*}
  \left(a+\frac{2}{m}\right)^2 + b^2 = \left(\frac{2}{m}\right)^2.
\end{equation*}
Thus for a fixed number $m$ of boundary lattice points the complex roots lie on the circle $|w+2/m|=2/m$. Moreover, infinitely many points on such a circle are zeros of   lattice polygons. For instance, for $m=3$ consider the triangle $T_q=\conv\{(-1,0)^\intercal, (1,-1)^\intercal, (0,q)^\intercal\}$. Then we have $\LE_2(T_q)=q+1/2$ and $\LE_1(T_q)=3/2$.

Real roots are just these with $g_1^2\geq 4g_2$ (cf.~\eqref{eq:tworoots}) or equivalently $(g_1/2-1)^2\geq g_2-g_1+1$. Hence by Pick's identity \eqref{eq:pick} and \eqref{eq:reciprocitylaw} we obtain ii). We remark that equality holds in ii) for the unit square and the triangle $\conv\{(0,0)^\intercal, (1,0)^\intercal, (0,2)^\intercal\}$ as well as  their integral multiplies.

For the third statement we use 
\cite[Theorem 2.2]{becketall:roots_ehrhart} where  it was shown that 
\begin{equation*}
  \Gamma(2)\subset \{-2,-1,-2/3\}\cup\left\{s=a+\imath b\in \C :
          -\frac{1}{2}\leq a <0, |b| \leq \sqrt{15}/6\right\}. 
\end{equation*} 
On the other hand we know by i) that all complex roots lie in the circle $|w+2/3|\leq 2/3$ and thus we get iii).

For iv) let $\gamma=p/q\in(0,1/2)$ with $p,q\in\N$ and let $\alpha=\gamma/(1-\gamma)=p/(q-p)$. For $k\in\N$ with $k\alpha\in\N$ we consider the $0$-symmetric hexagon   
\begin{equation*}
   H_k=\conv\{\pm(k,0)^\intercal, (\pm\,k\alpha,\pm 1)^\intercal\}.
\end{equation*}      
Writing $\lE_i$ instead of $\LE_i(H_k)$ we get $\lE_2=2(k\alpha+k)$ and $\lE_1=2(k\alpha+1)$. Hence for the roots of its Ehrhart polynomial we find 
\begin{equation*}
 -\gamma_{1,2}= -\frac{\lE_1}{2\,\lE_2}\left(1\pm\sqrt{1-4\frac{\lE_2}{\lE_1^2}}\right) = 
-\frac{1}{2}\frac{\alpha+1/k}{\alpha+1}\left(1\pm\sqrt{1-4\frac{k\,\alpha+k}{(k\,\alpha+1)^2}}\,\right).
\end{equation*}  
Thus for $k\to\infty$ the roots converge to $-\alpha/(\alpha+1)=-p/q$ and $0$, respectively. This shows that all points in $[-1/2,0]$ are cluster points of $\Gamma(2)$. It remains to show that there are no other cluster points. By iii) we know that any further cluster point has to be non real. So it suffices to show that for a given $\epsilon > 0$ there are only finitely many roots $\gamma=a+\imath\,b\in\Gamma(2)$ with $|b|>\epsilon$. 

Given an Ehrhart polynomial $\lE_2\,s^2+\lE_1\,s+1$ with  complex roots $-\gamma_{1,2}=a+\pm\imath\,b$ we  see by \eqref{eq:tworoots} that $|b|\leq 1/\sqrt{\lE_2}$. So a lower bound on $|b|$ gives an upper bound on $\lE_2$. Finally we observe, that there are only finitely many -- up to unimodular transformations -- different lattice polygons of prescribed volume.

If $P$ is a lattice polygons $P$ with $\LE(\inter(P))=0$ then by \eqref{eq:reciprocitylaw} we have $g_2-g_1+1=0$. Thus the roots of $P$ are given by $-1$ and $-2/(m-2)$,   where $m$ is the number of boundary lattice points (cf.~\eqref{eq:tworoots}). Hence we have verified v).

Next we show vi). 
In view of \eqref{eq:tworoots} the roots with real part $-1/2$ are given by the polygons with $\lE_2=\lE_1$ which by \eqref{eq:reciprocitylaw} is equivalent to $1=\LE(\inter(P))=g_1-g_2+1$. Thus, these are the lattice polygons with exactly one interior lattice point. Furthermore, by \eqref{eq:tworoots} we also have $4\,\lE_2-(\lE_1)^2\geq 0$. 
Setting again $\lE_1=m/2$ the last condition leads to $m\in\{3,\dots,8\}$. For $m=8$ we have only the real root $-1/2$ and hence, all together we have 11 roots. In fact there are lattice polygons with these roots. For $m=8$ we take the square of edge length $2$. Deleting successively the vertices of the square we obtain examples for $m=7,6,5,4$, the last one being the diamond. For $m=3$ we can take the triangle $T_1$ from part i) of this proof. 
By a result of Scott \cite{scott:latticepolygons} it is known that $g_2\leq 9/2$ for lattice polygons with only one interior lattice point and so we get $g_1\leq 9/2$. Hence there is only the case $m=9$ missing, i.e., the triangle $\conv\{(0,0)^\intercal, (3,0)^\intercal, (0,3)^\intercal\}$ with zeros $-2/3$, $-1/3$. 
\end{proof} 

\section{Remarks}
We conclude the paper with some open problems and questions. For a $0$-symmetric lattice polygon $P\in\mathcal{P}_0^2$ with roots $\gamma_i$ and successive minima $\lambda_i$ we know by \eqref{eq:minkroot} and \eqref{eq:main} that 
\begin{equation*}
 \frac{\lambda_1}{2}\frac{\lambda_2}{2}\leq \gamma_1\,\gamma_2 \text{ and } (\gamma_1+\gamma_2)^2 \leq \left(\frac{\lambda_1}{2}+\frac{\lambda_2}{2}\right)^2. 
\end{equation*} 
Thus  we get 
\begin{equation*}
 \gamma_1^2+\gamma_2^2 = (\gamma_1+\gamma_2)^2-2\gamma_1\gamma_2 \leq \left(\frac{\lambda_1}{2}+\frac{\lambda_2}{2}\right)^2 -2\frac{\lambda_1}{2}\frac{\lambda_2}{2}=\left(\frac{\lambda_1}{2}\right)^2+\left(\frac{\lambda_2}{2}\right)^2
\end{equation*}
and it is quite tempting to look for a generalisation to arbitrary dimensions, i.e., 
\begin{problem} Let $P\in\mathcal{P}_0^n$. Is it true that 
\begin{equation*}
  \sum_{i=1}^n \gamma_i(P)^2 \leq \sum_{i=1}^n \left(\frac{\lambda_i(P)}{2}\right)^2 \,\,?
\end{equation*}
\label{prob:one}
\end{problem}
In terms of the coefficients of the Ehrhart polynomial (cf.~\eqref{eq:coeffroots}) we may write 
\begin{equation*}
  \sum_{i=1}^n \gamma_i(P)^2 = \frac{\LE_{n-1}(P)^2-2\,\LE_n(P)\,\LE_{n-2}(P)}{\LE_n(P)^2}
\end{equation*} 
and so Problem \ref{prob:one} asks for a kind of Minkowski inequality  known in convexity (\cite[pp.~317]{schneider:93}) with respect to the functionals $\LE_i$. 

We also do not much about the size of the real roots of $0$-symmetric lattice polytopes. Numerical examples indicate that the lower bound $-1$ of Proposition \ref{prop:bound3} is not tight and can be replaced by $-\sqrt{3}/2$. Actually, we are not aware of $P\in\mathcal{P}_0^n$ with real roots less than -1. In the non-symmetric case it is known that a simplex provides the worst case. Thus  we would like to pose the problem
\begin{problem} Determine a sharp lower (and upper) bound on the real roots of $0$-symmetric lattice polytopes. 
\end{problem}
%-------------------------------------------------------------------------
\noindent 
{\it Acknowledgements.} We thank Iskander Aliev, Ulrich Betke and Jes\'us De Loera for helpful discussions. 
%---------------------------------bibliography-----------------------------
%\bibliographystyle{amsplain} 
%\bibliography{math_data}

\begin{thebibliography}{10}

\bibitem{bwz:65}
R.~P. Bambah, A.~C. Woods, and H.~Zassenhaus, \emph{Three proofs of
  {M}inkowski's second inequality in the geometry of numbers}, J. Austral.
  Math. Soc. \textbf{5} (1965), 453--462.

\bibitem{barvinok:ehrhart}
A.~Barvinok, \emph{Computing the {E}hrhart polynomial of a convex lattice
  polytope}, Discrete Comput. Geom. \textbf{12} (1994), no.~1, 35--48.

\bibitem{barvinok:polynomialtime}
\bysame, \emph{A polynomial time algorithm for counting integral points in
  polyhedra when the dimension is fixed}, Math.~Oper.~Res. \textbf{19} (1994),
  769--779.

\bibitem{barvinokpommersheim:algorithmic}
A.~Barvinok and J.E. Pommersheim, \emph{An algorithmic theory of lattice points
  in polyhedra}, New Perspectives in Algebraic Combinatorics (Berkeley, CA,
  1996-1997), Math. Sci. Res. Inst. Publ., vol.~38, 1994, pp.~91--147.

\bibitem{becketall:roots_ehrhart}
M.~Beck, J.~{De Loera}, M.~Develin, J.~Pfeifle, and R.P. Stanley,
  \emph{Coefficients and roots of {E}hrhart polynomials}, Contemp. Math.
  \textbf{374} (2005), 15--36.

\bibitem{beckpixton:birkhoff}
M.~Beck and D.~Pixton, \emph{The {E}hrhart polynomial of the {B}irkhoff
  polytope}, Discrete Comput.~Geom. \textbf{30} (2003), no.~4, 623--637.

\bibitem{beckrobins:discrete_continously}
M.~Beck and S.~Robins, \emph{Computing the continuous discretely: Integer-point
  enumeration in polyhedra}, Springer, (to appear), Preprint at
  \url{http://math.sfsu.edu/beck/papers/ccd.html}.

\bibitem{betkegritzmann:valuations}
U.~Betke and P.~Gritzmann, \emph{An application of valuation theory to two
  problems of discrete geometry}, Discrete Math. \textbf{58} (1986), 81--85.

\bibitem{bhw:93}
U.~Betke, M.~Henk, and J.~M. Wills, \emph{Successive-minima-type inequalities},
  Discrete Comput. Geom. \textbf{9} (1993), no.~2, 165--175.

\bibitem{betkekneser:latticehadwiger}
U.~Betke and M.~Kneser, \emph{Zerlegungen und {B}ewertungen von
  {G}itterpolytopen}, J.~{R}eine {A}ngew.~{M}ath. \textbf{358} (1985),
  202--208.

\bibitem{betke-mcMullen:85}
U.~Betke and P.~McMullen, \emph{Lattice points in lattice polytopes}, Monatsh.
  Math. \textbf{99} (1985), no.~4, 253--265.

\bibitem{bump:crosspolytope}
D.~Bump, K.-K. Choi, P.~Kurlberg, and J.~Vaaler, \emph{A local {R}iemann
  hypothesis {I}}, Math. Z. \textbf{233} (2000), no.~1, 1--19.

\bibitem{davenport:39}
H.~Davenport, \emph{Minkowski's inequality for the minima associated with a
  convex body}, Quarterly J.~Math. \textbf{10} (1939), 119--121.

\bibitem{latte:program}
J.~{De Loera}, D.~Haws, R.~Hemmecke, and P.~Huggins, \emph{A user's guide for
  \url{latte} v1.1, software package \url{latte}}, 2004, available at
  \url{http://www.math.ucdavis.edu/~latte}.

\bibitem{deloera:effecticve_counting}
J.~{De Loera}, R.~Hemmecke, J.~Tauzer, and R.~Yoshida, \emph{Effective lattice
  point counting in rational convex polytopes}, J.~Symb.~Comput. \textbf{38}
  (2004), no.~4, 1273--1302.

\bibitem{diazrobins:latticepolytope}
R.~Diaz and S.~Robins, \emph{The {E}hrhart polynomial of a lattice polytope},
  Ann.~of Math. \textbf{145} (1997), no.~3, 503--518, Erratum in 146:1 (1997),
  237.

\bibitem{ehrhart:polynomial}
E.~Ehrhart, \emph{Sur les poly\`edres rationnels homoth\'etiques \`a n
  dimensions}, C. R. Acad. Sci., Paris, S\'er. A \textbf{254} (1962), 616--618.

\bibitem{ehrhart:polynomialII}
\bysame, \emph{Sur un probl\`eme de g\'eom\'etrie diophantienne lin\'eaire},
  J.~{R}eine {A}ngew.~Math. \textbf{227} (1967), 25--49.

\bibitem{ehrhart:reciprocity}
\bysame, \emph{Sur la loi de r\'eciprocit\'e des poly\`edres rationnels}, C. R.
  Acad. Sci., Paris, S\'er. A \textbf{266} (1968), 695--697.

\bibitem{GriWil:latticepoints}
P.~Gritzmann and J.M. Wills, \emph{Lattice points}, Handbook of convex geometry
  (P.M. Gruber and J.M. Wills, eds.), vol.~B, North-Holland, Amsterdam, 1993.

\bibitem{gruber-lekkerkerker:87}
P.~M. Gruber and C.~G. Lekkerkerker, \emph{Geometry of numbers}, second ed.,
  vol.~37, North-Holland Publishing Co., Amsterdam, 1987.

\bibitem{gruenbaum:polytopes}
B.~Gr\"unbaum, \emph{Convex polytopes}, 2nd ed., Springer, New-York, 2003,
  Second edition prepared by V.~Kaibel, V.~Klee and G.~M.~Ziegler.

\bibitem{henk:90}
M.~Henk, \emph{Inequalities between successive minima and intrinsic volumes of
  a convex body}, Monatsh. Math. \textbf{110} (1990), no.~3-4, 279--282.

\bibitem{henk:02}
\bysame, \emph{Successive minima and lattice points}, Rend. Circ. Mat. Palermo
  (2) Suppl. (2002), no.~70, part I, 377--384.

\bibitem{hibi:algebraic_combinatorics}
T.~Hibi, \emph{Algebraic combinatorics on convex polytopes}, Carslaw
  Publications, Glebe, Australia, 1992.

\bibitem{hibi:92}
T.~Hibi, \emph{Dual polytopes of rational convex polytopes}, Combinatorica
  \textbf{12} (1992), no.~2, 237--240.

\bibitem{lagarias:pointlattices}
J.C. Lagarias, \emph{Point lattices}, Handbook of Combinatorics (R.L. Graham,
  M.~Gr\"otschel, and L.~Lov\'asz, eds.), vol.~A, North-Holland, Amsterdam,
  1995.

\bibitem{lagariasziegler:volume}
J.C. Lagarias and G.M. Ziegler, \emph{Bounds for lattice polytopes containing a
  fixed number of interior points in a sublattice}, Canad.~J.~Math. \textbf{43}
  (1991), 1022--1035.

\bibitem{liu:cyclic_ehrhart}
F.~Liu, \emph{Ehrhart polynomials of cyclic polytopes}, J.~Comb.~Theory, Ser.~A
  \textbf{111} (2005), 111--127.

\bibitem{minkowski:96}
H.~Minkowski, \emph{Geometrie der {Z}ahlen}, Teubner, Leipzig-Berlin, 1896,
  Reprinted: Johnson, New York, 1968.

\bibitem{mordell:dedkind}
L.J. Mordell, \emph{Lattice points in tetrahedron and generalized {D}edekind
  sums}, J.~Indian Math.~Soc. (N.S.) \textbf{15} (1951), 41--46.

\bibitem{mustatepayne:bettyehrhart}
M.~Mustata and S.~Payne, \emph{Ehrhart polynomials and stringy {B}etti
  numbers}, \url{http://arxiv.org/abs/math.AG/0505054}.

\bibitem{pick:lattice}
G.A. Pick, \emph{Geometrisches zur {Z}ahlenlehre}, Sitzungsber.~Lotus Prag
  \textbf{19} (1899), 311--319.

\bibitem{pommersheim:tetrahedra}
J.E. Pommersheim, \emph{Toric varieties, lattice points and {D}edekind sums},
  Math.~Ann. \textbf{295} (1993), no.~1, 1--24.

\bibitem{reeve:simplex}
J.E. Reeve, \emph{On the volume of lattice polyhedra}, Proc.~London Math.~Soc.
  \textbf{7} (1957), no.~3, 378--395.

\bibitem{rodriguez:zeros}
F.~Rodriguez-Villegas, \emph{On the zeros of certain polynomials}, Proc. Amer.
  Math. Soc. \textbf{130} (2002), 2251--2254.

\bibitem{schneider:93}
R.~Schneider, \emph{Convex bodies: the {B}runn-{M}inkowski theory},
  Encyclopedia of Mathematics and its Applications, vol.~44, Cambridge
  University Press, Cambridge, 1993.

\bibitem{scott:latticepolygons}
P.R. Scott, \emph{On convex lattice polygons}, Bull.~Austral.~Math.~Soc.
  \textbf{15} (1976), 395--399.

\bibitem{siegel:89}
C.~L. Siegel, \emph{Lectures on the geometry of numbers}, Springer-Verlag,
  Berlin, 1989.

\bibitem{stanley:nonnegative}
R.~P. Stanley, \emph{Decompositions of rational convex polytopes},
  Ann.~Discrete Math. \textbf{6} (1980), 333--342.

\bibitem{stanley:enumcombI}
\bysame, \emph{Enumerative combinatorics. vol. 1}, Cambridge {S}tudies in
  {A}dvanced {M}athematics, vol.~49, Cambridge University Press, Cambridge,
  1997, Corrected reprint of the 1986 original.

\bibitem{weyl:42}
H.~Weyl, \emph{On geometry of numbers}, Proc. London Math. Soc. (2) \textbf{47}
  (1942), 268--289.

\bibitem{wills:90a}
J.~M. Wills, \emph{Kugellagerungen und {K}onvexgeometrie}, Jahresber. Deutsch.
  Math.-Verein. \textbf{92} (1990), no.~1, 21--46.

\bibitem{wills:90c}
\bysame, \emph{Minkowski's successive minima and the zeros of a
  convexity-function}, Monatsh. Math. \textbf{109} (1990), no.~2, 157--164.

\bibitem{wills:minkowskianalog}
J.M. Wills, \emph{On an analog to {M}inkowski's lattice point theorem}, The
  Geometric Vein: The {C}oxeter {F}estschrift (C.~Davis, B.~Gr\"unbaum, and
  F.A. Sherk, eds.), Springer, New York, 1982, pp.~285--288.

\end{thebibliography}
\def\cprime{$'$} \def\cprime{$'$} \def\cprime{$'$} \def\cprime{$'$}
  \def\cprime{$'$} \def\cprime{$'$} \def\cprime{$'$}
\providecommand{\bysame}{\leavevmode\hbox to3em{\hrulefill}\thinspace}
\providecommand{\MR}{\relax\ifhmode\unskip\space\fi MR }
% \MRhref is called by the amsart/book/proc definition of \MR.
\providecommand{\MRhref}[2]{%
  \href{http://www.ams.org/mathscinet-getitem?mr=#1}{#2}
}
\providecommand{\href}[2]{#2}

%--------------------------------------------------------------------------
\end{document}